\documentclass{amsart}
\usepackage{fourier} 
\usepackage{latexsym}
\usepackage{amssymb}
\usepackage{amsmath}
\usepackage{amscd}
\usepackage{float}
\usepackage{mathrsfs} 
\usepackage{enumerate}
\usepackage{tikz}

\tikzset{node distance=1.5cm, auto}

\newtheorem{theorem}{Theorem}
\newtheorem{proposition}{Proposition}
\newtheorem{lemma}{Lemma}
\newtheorem{corollary}{Corollary}

\theoremstyle{definition}
\newtheorem{definition}{Definition}
\theoremstyle{remark}

\def\<{\langle}
\def\>{\rangle}

\def\RR{\mathbb{R}}

\newcommand{\sgn}{\operatorname{sgn}}

\title{Tropical Multiplier Sequences}
\date{\today}
\author{Jens Forsg{\aa}rd}
\address{Department of Mathematics \\ Stockholm University \\
SE-106 91 Stockholm, Sweden.}
\email{jensf@math.su.se}

\begin{document}

\begin{abstract}
We prove that the diagonal operator defined by a 
positive sequence preserves tropical and central indices
if and only if the sequence is log-concave. 
In particular we obtain an
elementary proof of that such an operator
preserves the set of sign-independently real-rooted
polynomials if and only if the sequence
is log-concave.
\end{abstract}

\maketitle

\vspace{-1cm}
\renewcommand{\abstractname}{R{\'e}sum{\'e}}
\begin{abstract}
Nous d{\'e}montrons que l'operateur diagonal defini par une s{\'e}quence 
positive pr{\'e}serve les index tropicales et centrales si et seulement 
si la s{\'e}quence est log-concave. En particulier nous obtenons une 
d{\'e}montration {\'e}lementaire du fait qu'un tel op{\'e}rateur pr{\'e}serve l'ensemble 
de p{\^o}lynomes {\`a} racines r{\'e}elles, ind{\'e}pendant du signe, si et seulement si 
la s{\'e}quence est log-concave.
\end{abstract}
\vspace{0.5cm}

\section{Results}

Consider the polynomial ring $\RR[z]$ consisting of all real 
univariate polynomials.
Let $\gamma = \{\gamma_n\}_{n=0}^\infty$ be a sequence of real
numbers, to which we associate the diagonal operator 
$T_\gamma\colon \RR[z]\rightarrow \RR[z]$ 
defined by $z^n \mapsto \gamma_n z^n$, for $n=0, 1, \dots$,
and extended to $\RR[z]$ by linearity.
The sequence $\gamma$ is said to be a multiplier sequence 
if $T_\gamma$ preserves the set
of real-rooted polynomials, see \cite{CC04}
for background information and applications.

We borrow the following notation from Wiman--Valiron theory,
see, e.g., \cite{Hay74}.
A non-negative integer $m$ is said to be a \emph{central index} of the 
polynomial
\[
f(z) = \sum_{n=0}^d a_n z^n 
\]
if there exists a number $z_m\geq 0$ such that
\begin{equation}
\label{eqn:Lopsidedness}
 |a_m|z_m^m \geq \sum_{n\neq m}|a_n|z_m^n.
\end{equation}
Condition \eqref{eqn:Lopsidedness} has recently appeared in the context of
amoebas, see, e.g., \cite{Rul03}.
To relate \eqref{eqn:Lopsidedness} to real-rootedness, we recall 
that a real polynomial is said to be sign-independently real-rooted if any
polynomial obtained by arbitrary sign changess of its coefficients
is real-rooted, see \cite{PRS11}. 

\begin{proposition}
\label{pro:SignIndependently}
 A real polynomial $f$ of degree $d$ is sign-independently real-rooted
 if and only if each index $n=0, \dots, d$ is a central index of $f$.
\end{proposition}

In this article, we discuss sequences $\gamma$ that preserve
the set of 
central indices. In addition, we introduce the following notion.
A non-negative integer $m$ is said to be a \emph{tropical index} of $f$ if there exists 
a number $z_m\geq 0$ such that
\begin{equation}
\label{eqn:Tropical}
 |a_m|z_m^m \geq \max_{n\neq m}\,|a_n|z_m^n.
\end{equation}
Notice that \eqref{eqn:Tropical} is the analogue of 
\eqref{eqn:Lopsidedness} when the right
hand side of \eqref{eqn:Lopsidedness} is 
interpreted as a tropical sum.
A polynomial $f$ of degree $d$
is said to be \emph{tropically real-rooted}
if and only if each index $n=0, \dots, d$ is a tropical index
of $f$.

As the definition of central and tropical indices only depend
on the moduli $|a_n|$,
they are immediately extended to complex polynomials.
However, for simplicity, we will henceforth assume that $a_n \geq 0$
for all $n$, and we will consider only positive sequences $\gamma$.
Such a sequence is said to be \emph{log-concave} if $\gamma_n^2 \geq \gamma_{n-1}\gamma_{n+1}$ 
for all $n$.
In \cite{PRS11} it was proven, using discriminant amoebas,
that the diagonal operator $T_\gamma$ associated to the sequence 
$\gamma$ preserves the set of sign-independently
real-rooted polynomials if and only if $\gamma$ is log-concave.
For this reason, log-concave sequences are said to be
\emph{multiplier sequences of the third kind}
(we will say that $\gamma$ is a \emph{tropical multiplier sequence}). 

\begin{definition}
 A sequence $\gamma$ is said to be a tropical (resp.\ central) index preserver
 if for each polynomial $f$ the set of tropical (resp.\ central) indices
 of $f$ is a subset of the set of tropical (resp.\ central) indices
 of $T_\gamma[f]$.
\end{definition}

Our main results are as follows.

\begin{theorem}
\label{thm:TropicalIndexPreserver}
A positive sequence $\gamma$ is a tropical index preserver 
if and only if it is log-concave.
\end{theorem}

\begin{theorem}
\label{thm:CentralIndexPreserver}
A positive sequence $\gamma$ is a central index preserver if and only if 
it is log-concave.
\end{theorem}

As a corollary, we obtain
an elementary proof of \cite[Theorem 1]{PRS11}
as requested in \cite[Problem  2]{PRS11}.

\begin{corollary}
\label{cor:SignIndependently}
A positive sequence $\gamma$ preserves the set of sign-independently 
real-rooted polynomials if and only if it is log-concave.
\end{corollary}

Finally, let us state our main lemma.

\begin{lemma}
\label{lem:LogConcaveTropical}
A positive sequence $\gamma$ is log-concave if and only if for each $d$ the 
polynomial 
\[
P_\gamma(z) = \sum_{n=0}^d \gamma_n z^n
\]
is tropically real-rooted.
\end{lemma}

Using Lemma \ref{lem:LogConcaveTropical}, we could rephrase
Theorems \ref{thm:TropicalIndexPreserver} and \ref{thm:CentralIndexPreserver}
in a manner similar to the classical result of 
P{\'o}lya and Schur \cite{PS14}.
Namely, a positive sequence $\gamma$
is a central and tropical index preserver if and
only if the \emph{tropical symbol} $P_\gamma(z)$ is
tropically real-rooted.

\section{Acknowledgements}
I would like to thank B.\ Shapiro for comments and encouragement, 
and C.\ Esp{\'i}ndola for the translation into French.

\section{Proofs}

\begin{proof}[Proof\/ of\/ Lemma \ref{lem:LogConcaveTropical}]
Assume first that $\gamma$ is log-concave. For each $m \geq 1$ define $z_m$ by
$z_m = \sqrt{\gamma_{m-1}/\gamma_{m+1}}$.
Then,
\[
\frac{z_{m+1}}{z_m} 
= \frac{\gamma_{m}}{\sqrt{\gamma_{m-1}\gamma_{m+1}}}\frac{\gamma_{m+1}}{\sqrt{\gamma_m\gamma_{m+2}}} \geq 1,
\]
so that $\{z_m\}_{m=1}^\infty$ is a non-decreasing sequence of positive real
numbers. Further more,
\[
\frac{\gamma_m z_m^m}{\gamma_{m-1} z_m^{m-1}} = 
\frac{\gamma_m z_m^m}{\gamma_{m+1} z_m^{m+1}} = 
\frac{\gamma_m}{\sqrt{\gamma_{m-1}\gamma_{m+1}}} \geq 1.
\]
Since both binomials $\gamma_n z^n - \gamma_{n+1}z^{n+1}$ and 
$\gamma_n z^n - \gamma_{n-1}z^{n-1}$ have exactly 
one positive real root, we conclude that
$\gamma_n z_m^n \geq \gamma_{n+1}z_m^{n+1}$ if $n\geq m$ and that
$\gamma_n z_m^n \geq \gamma_{n-1}z_m^{n-1}$ if $n\leq m$.
Hence,
\[
\gamma_m z_m^m \geq \max_{n\neq m}\, \gamma_n z_m^n.
\]

For the converse, assume that $\gamma$ is not log-concave. That is, there exists
an index $m$ for which $\gamma_m^2 < \gamma_{m-1}\gamma_{m+1}$.
Then, for $z\geq 0$,
\[
\gamma_m z^m < \sqrt{\gamma_{m-1}z^{m-1}\, \gamma_{m+1}z^{m+1}}
\leq \max \left(\gamma_{m-1}z^{m-1}, \gamma_{m+1}z^{m+1}\right).
\]
In particular, $m$ is not a tropical index of $P_\gamma(z)$.
\end{proof}

\begin{proof}[Proof\/ of\/ Theorem \ref{thm:TropicalIndexPreserver}]
Assume first that $\gamma$ is log-concave.
Let $m$ be a tropical index of $f$, and let $z_m\geq 0$ be such that
\[
a_m z_m^m \geq \max_{n\neq m}\, a_n z_m^n.
\]
By Lemma \ref{lem:LogConcaveTropical} we can find a $\zeta_m$ such that
\[
\gamma_m\zeta_m^m \geq \max_{n\neq m}\, \gamma_n\zeta_m^n.
\]
Then
\[
\gamma_ma_m(z_m \zeta_m)^m =\gamma_mz_m^m\,a_m \zeta_m^m \geq 
\gamma_nz_m^n\,a_n \zeta_m^n
\]
for all $n$. Hence, $m$ is a tropical index of $T_\gamma[f]$.

For the converse, it suffices to consider the polynomials $1+z+\dots + z^d$,
which is tropically real-rooted for all $d$, and 
use Lemma \ref{lem:LogConcaveTropical}.
\end{proof}

\begin{proof}[Proof\/ of\/ Theorem \ref{thm:CentralIndexPreserver}]
Assume first that $\gamma$ is log-concave, and
let $\zeta_m$ be as in the proof of Theorem \ref{thm:TropicalIndexPreserver}.
Let $m$ be a central index of $f$, and let $z_m$ be such that
\[
a_mz_m^m \geq \sum_{n\neq m} a_n z_m^n.
\]
Then,
\[
\gamma_m a_m (z_m\zeta_m)^m \geq \sum_{n\neq m}\gamma_m \zeta_m^m a_n z_m^n
\geq \sum_{n\neq m}\gamma_n \zeta_m^n a_n z_m^n,
\]
implying that $m$ is a central index of $T_\gamma[f]$.

For the converse, assume that $\gamma_m^2 < \gamma_{m-1}\gamma_{m+1}$,
and consider the action of $T_\gamma$ on the trinomial $z^{m-1} + 2 z^m + z^{m+1}$.
\end{proof}

\begin{proof}[Proof\/ of\/ Proposition \ref{pro:SignIndependently}]
To prove the \emph{only if}-part, it suffices to
choose $z_m$ as the mean of the two positive roots of the polynomial
\[
|a_m|z^m - \sum_{n\neq m} |a_n| z^n,
\]
which exists by assumption.
For the \emph{if}-part, choose arbitrary signs of the coefficients of $f$.
We note that the criterion \eqref{eqn:Lopsidedness} 
implies that
\[
\sgn(f(z_m)) = \sgn(a_mz_m^m) = \sgn(a_m),
\]
for $z>0$. Using additionally Descartes' rule of signs, we conclude that
the number of positive roots of $f$ is equal to the number of sign
changes in the sequence $\{a_n\}_{n=0}^d$. Similarly, the number 
of negative roots of $f$ is equal to the number of sign changes in
the sequence $\{a_n(-1)^n\}_{n=0}^d$. 
As $a_n\neq 0$ for each $n$, these two numbers
sums up to $d$, implying that $f(z)$ is real-rooted.
Since the signs of the coefficients were chosen arbitrary,
we are done.
\end{proof}

\begin{proof}[Proof\/ of\/ Corollary \ref{cor:SignIndependently}]
It follows from Proposition \ref{pro:SignIndependently} that a positive sequence preserves the set of sign-independently
real-rooted polynomials if and only if it preserves central indices, and
it follows from Theorem \ref{thm:CentralIndexPreserver} that a positive sequence preserves central indices
if and only if it is log-concave.
\end{proof}


\begin{thebibliography}{99}
\bibliographystyle{amsplain}

\bibitem{CC04}
Craven, T.\ and Csordas, G.,
\emph{Composition theorems, multiplier sequences and complex zero decreasing sequences},
Value distribution and related topics,
Adv.\ Complex Anal.\ Appl., vol.\ 3,
Kluwer Acd., Publ., Boston, MA, 2004, pp.\ 131--166.

\bibitem{Hay74}
Hayman, W.\ K.,
\emph{The local growth of power series: a survey of the Wiman--Valiron method},
Canad.\ Math.\ Bull.\ \textbf{17} (1974), no.\ 3, 317--358.

\bibitem{PRS11}
Passare, M., Rojas, J.\ M., and Shapiro, B.,
\emph{New multiplier sequences via discriminant amoebae},
Mosc.\ Math.\ J.\ \textbf{11} (2011), no.\ 3, 547--560, 631.

\bibitem{PS14}
P{\'o}lya, G.\ and Schurm J., 
\emph{{\"U}ber zwei Arten von Fasterfolgen},
J.\ Reine Angew.\ Math.\ \textbf{144} (1914),
89--133.

\bibitem{Rul03}
Rullg{\aa}rd, H.,
\emph{Topics in geometry, analysis, and inverse problems},
Doctoral thesis, Stockholms universitet, Stockholm, 2003.

\end{thebibliography}
\end{document}